\documentclass[leqno,11pt]{amsart}
\usepackage{amsmath,amsthm, amscd, amssymb, amsfonts}
\usepackage[all]{xy}
\usepackage{hhline}
\usepackage[dvips, dvipsnames, usenames]{color}

\newcommand{\Ind}{\operatorname{Ind}}

\newcommand{\oct}{\mathfrak O}

\newcommand\oc{{\mathcal O }}
\newcommand\ocs{{\mathcal O_{\sigma} }}
\newcommand\toba{{\mathfrak B }}

\newcommand{\trid}{\triangleright}

\newcommand{\ku}{\mathbb C}

\newcommand{\Z}{{\mathbb Z}}

\newcommand{\Q}{{\mathcal Q}}

\newcommand{\D}{{\mathcal D}}

\newcommand{\B}{{\mathcal B}}

\newcommand{\Oc}{{\mathcal O}}

\newcommand\ot{\otimes}

\newcommand\sgn{\operatorname{sgn}}

\theoremstyle{plain}

\newtheorem{maintheorem}{Theorem}
\newtheorem{lema}{Lemma}[section]

\newtheorem{prop}[lema]{Proposition}
\newtheorem{claim}{Claim}

\theoremstyle{definition}

\theoremstyle{remark}
\newtheorem{obs}[lema]{Remark}

\newcommand\id{\operatorname{id}}

\newcommand\Id{\operatorname{Id}}

\newcommand\s{\mathbb S}

\newcommand\sm{\mathbb S_m}

\newcommand\sigmae{\sigma_{e}}
\newcommand\sigmao{\sigma_{o}}
\newcommand\so{\sigmao}

\def\pf{\begin{proof}}
\def\epf{\end{proof}}

\theoremstyle{remark}

\begin{document}

\renewcommand{\baselinestretch}{1.2}

\thispagestyle{empty}

\title[pointed Hopf algebras: symmetric groups]{On pointed Hopf algebras associated \\ with the symmetric groups}

\author[andruskiewitsch, fantino and zhang]{Nicol\'{a}s Andruskiewitsch, Fernando Fantino \\ and Shouchuan Zhang}
\address{\noindent Facultad de Matem\'{a}tica, Astronom\'{i}a y F\'{i}sica\\
Universidad Nacional de C\'{o}rdoba \\ CIEM - CONICET,
(5000) Ciudad Universitaria \\
C\'{o}rdoba \\Argentina}\email{andrus@famaf.unc.edu.ar}
\email{fantino@famaf.unc.edu.ar}
\address{\noindent Department  of Mathematics, Human University\\ Changsha
410082, \ P.R. China } \email{z9491@yahoo.com.cn}
\thanks{This work was partially
supported by CONICET, ANPCyT and Secyt (UNC)}

\subjclass[2000]{16W30; 17B37}
\date{\today}

\begin{abstract}
It is an important open problem whether the dimension of the
Nichols algebra $\toba(\oc,\rho)$ is finite when $\oc$ is the
class of the transpositions and $\rho$ is the sign representation,
with $m\geq 6$. In the present paper, we discard most of the other
conjugacy classes showing that very few pairs $(\oc,\rho)$ might
give rise to finite-dimensional Nichols algebras.
\end{abstract}

\maketitle

\section{Introduction and Main Result}\label {0}

\subsection{The context}

This paper contributes to the classification of finite-dimensional
pointed \emph{complex} Hopf algebras $H$ whose group of
group-likes $G(H)$ is isomorphic to $\sm$.  Suppose we want to
classify finite-dimensio\-nal pointed Hopf algebras with a fixed
$G(H)=G$. As explained in \cite{AS3}, the crucial step is to
determine when the Nichols algebra of a Yetter-Drinfeld module
over $G$ is finite-dimensional or not.  Recall that irreducible
Yetter-Drinfeld modules over $G$ are determined by a conjugacy
class $\oc$ of $G$ and an irreducible representation of the
centralizer $G^{\sigma}$ of a fixed $\sigma \in \oc$. Let
$M(\oc,\rho)$ be the corresponding Yetter-Drinfeld module and let
$\toba(\oc,\rho)$ denote its Nichols algebra.

\subsection{Statement of the main result}

Assume that $G=\sm$, $3\leq m \in \mathbb N$. We fix $\sigma \in
\mathbb S_m$ of type $(1^{n_1},2^{n_2},\dots,m^{n_m})$, which
means that in the decomposition of $\sigma$ as product of disjoint
cycles appear $n_j$ cycles of length $j$, for every $j$, $1\leq j
\leq m$. We will write
\begin{align}\label{eqn:descomp}
\sigma=A_1 \cdots A_m,
\end{align}
where $A_j=A_{1,j} \cdots A_{n_j,j}$ is the product of the $n_j
>0$ disjoint $j$-cycles $A_{1,j}$, \dots, $A_{n_j,j}$ of
$\sigma$. We omit $A_j$ when $n_j =0$. The even and the odd parts
of $\sigma$ are
\begin{align}\label{sigma:even:odd}
\sigmae:=\prod_{j \text{ even}}A_j, \qquad \sigma_{o}:=\prod_{1 <
j \text{ odd}}A_j.
\end{align}
Thus, $\sigma = A_1\sigmae \sigma_{o}$. By abuse of notation we
shall say that $\sigma$ has type $(1^{n_1}, 2^{n_2}, \dots, \so)$
to mean that the type of the $\so$ entering in $\sigma$ is
arbitrary. Here is our main Theorem. See below for unexplained
notation; in particular, see \eqref{formrho2} for the meaning of
$\rho_j$.

\begin{maintheorem}\label{maintheorem}
Let $\sigma\in \mathbb S_m$ be of type $(1^{n_1}, 2^{n_2}, \dots,
m^{n_m})$, let $\Oc$ be the conjugacy class of $\sigma$ and let
$\rho=(\rho,V) \in \widehat{\mathbb S_m^{\sigma}}$. Assume that
$\dim\toba(\Oc, \rho) < \infty$. Then $q_{\sigma\sigma}=-1$ and
some of the following hold:

\renewcommand{\theenumi}{\roman{enumi}}   \renewcommand{\labelenumi}{(\theenumi)}

\begin{enumerate}
\item $(1^{n_1}, 2)$, $\rho_1 =\sgn$ or $\epsilon$, $\rho_2 =\sgn$.
\item $(2, \so)$, $\so\neq \id$, $\rho_2 =\sgn$, $\rho_j=\Ind\chi_{(0,\dots,0)}\otimes \mu_j$, for all $j>1$ odd.
\item $(1^{n_1}, 2^3)$, $\rho_1 =\sgn$ or $\epsilon$, $\rho_2 = \chi_{(3)}\otimes \epsilon$ or $\chi_{(3)}\otimes \sgn$.\\
Furthermore, if $n_1>0$, then $\rho_2 = \chi_{(3)}\otimes \sgn$.
\item $( 2^5)$, $\rho_2 = \chi_{(5)}\otimes \epsilon$ or $\chi_{(5)}\otimes \sgn$.
\item $(1^{n_1}, 4)$, $\rho_1 =\sgn$ or $\epsilon$, $\rho_4=\chi_{(-1)}$.
\item $(1^{n_1}, 4^2)$, $\rho_1 =\sgn$ or $\epsilon$, $\rho_4=\chi_{(i,i)}\otimes \sgn$ or
$\chi_{(-i,-i)}\otimes \sgn$.
\item $(2,4)$, $\rho=\sgn \otimes \epsilon$ or
$\rho=\epsilon \otimes \chi_{(-1)}$.
\item $(2,4^2)$, $\rho_2=\epsilon$, $\rho_4=\chi_{(i,i)}\otimes \sgn$ or
$\chi_{(-i,-i)}\otimes \sgn$.
\item $(2^{2},4)$, $\deg\rho_2=1$, $\rho_4=\chi_{(-1)}$.
\end{enumerate}
\end{maintheorem}

We stress that in most of the cases en the previous statement,
whether the dimension of the corresponding Nichols algebra is
finite is an open problem; the point of the Theorem is to discard
the cases not comprised in (i) to (ix).

\subsection{Proof of the main result}

In the previous papers \cite{AF1, AZ} and in a preliminar version
of the present paper, the idea was to look at abelian subracks.
But we have found that the techniques based on non-abelian
subracks presented in \cite{AF3}-- consequences of the results in
\cite{AHS}-- conduct faster to a more complete analysis of Nichols
algebras over symmetric groups. We use also here another
technique, an extension of \cite{G1}, consisting in finding a
suitable braided vector subspace of diagonal type not supported by
an abelian subrack but ``transversal" -- see Proposition
\ref{te:2^k}. We now outline the proof of the main Theorem,
addressing to previous papers or to results in the present paper
for proofs of the different steps.

\

\noindent\emph{Proof.} We state the different arguments that
reduce the class of possible Nichols algebras with finite
dimension.

\renewcommand{\theenumi}{\alph{enumi}}   \renewcommand{\labelenumi}{(\theenumi)}
\begin{enumerate}
\item $q_{\sigma\sigma}=-1$ and $\sigma$ has even order by \cite[2.2]{AZ}, cf. Lemma \ref{le:odd}
below.

\medbreak\item If $j\geq 6$ is even, then $n_j=0$; this follows
from \cite[Ex. 2.10]{AF3} for $j$ having an odd divisor, and from
Proposition \ref{te:2^k} for $j$ a power of 2. Hence $\deg \rho_1
= 1$ by Propositions \ref{prop:degrho1>1:n2} and
\ref{prop:degrho1>1:n4}; that is, $\rho_1 =\sgn$ or $\epsilon$.

\medbreak\item $n_{4} \leq 2$, by \cite[Ex. 3.10]{AF3}; $n_2 \leq
5$, by \cite[Ex. 3.13]{AF3}.

\medbreak\item $\deg\rho_2=\deg\rho_4=1$, by Proposition 3.4 and
\cite[Prop. 2.6]{AZ}.

\medbreak\item Assume  that  \emph{there exists $j\geq 3$ such
that $n_j
> 0$.}
Then $n_2\leq 2$ by \cite[Ex. 3.12]{AF3}. If $n_2>0$, then $n_1 =
0$ by \cite[Ex. 3.9]{AF3}. Moreover, $(2^2, 4^2, \so)$ is excluded
by Proposition \ref{prop:n2=n4=2}.

\medbreak\item If $n_4 >0$, then $\so$ is trivial by Prop.
\ref{prop:n4>0yod}. The restrictions on the representations follow
from Proposition \ref{prop:n4=2}.

\medbreak\item Assume next that  \emph{for any $j\geq 3$, $n_j =
0$}, so that $\sigma$ is of type $(1^{n_1}, 2^{n_2})$. Then
$n_2\neq 2$ by \cite[Th. 2.7]{AZ} and $n_2\neq 4$ by \cite[Th. 1
(B) (i)]{AF1} together with \cite[Prop. 2.6]{AZ}. Also, $\rho_2$
should be $\chi_{(n_2)}\otimes \epsilon$ or $\chi_{(n_2)}\otimes
\sgn$ by \cite[Th. 1 (B) (ii)]{AF1}. Now, by Lemma \ref{prop:n1>0}
we have the restrictions on the characters in (iii) and (iv). On
the other hand, the claim in (ii) for the representations
$\rho_j$, $j>1$ odd, follows from Lemma \ref{lema:odd:gral}.\qed
\end{enumerate}

\subsection{Comments on the cases left open} Let us say that $M(\oc,\rho)$ has \emph{negative braiding} if the
Nichols algebra of any braided subspace corresponding to an
abelian subrack is (twist-equivalent to) an exterior algebra.

\begin{enumerate}
\item[(I)] In the cases (i)-(ix) of the main Theorem, there is no family of
type $\oct^{(2)}$ nor $\D_p^{(2)}$, for any odd prime $p$, inside
the respective conjugacy classes. Besides, there is no transversal
subrack of type $\D_4^{(2)}$. Moreover, it is easy to see that the
corresponding braiding is negative in all cases.

\item[(II)] We know that $\dim \toba(\oc,\rho)< \infty$ in the
following cases: $\oc$ is the class of the transpositions, $\rho =
\id\ot\sgn$ or $\epsilon\ot\sgn$ (see the notation below) and
$m\leq 5$  \cite{ms, FK, G2}; or $\oc$ is the class of the
4-cycles, $\rho = \chi_{(-1)}$ and $m = 4$  \cite[Th. 6.12]{AG2}.
It is an important open problem whether the dimension of
$\toba(\oc,\rho)$ is finite when $\oc$ is the class of the
transpositions and $m\geq 6$.

\item[(III)] The Yetter-Drinfeld modules over the group algebra $\ku
G$, with $G$ a finite group, are semisimple. The Nichols algebra
of a finite-dimensional reducible Yetter-Drinfeld module over
$\s_m$ is always in\-fi\-ni\-te-dimensional -- see \cite{HS}, see
also \cite[Section 4]{AHS}.

\end{enumerate}

\section{Preliminaries}\label{conventions}

\subsection{Generalities}
We  follow the conventions in \cite{AZ,AF1}. We denote by
$\widehat{G}$ the set of isomorphism classes of irreducible
representations of a finite group $G$. We  use the rack notation
$x\trid y:=xyx^{-1}$. We set $\omega _n : = e ^ {\frac {2\pi i}
{n}}$, where $i= \sqrt{-1}$.

Let $\sigma \in G$, $\Oc_{\sigma}$ the conjugacy class of $\sigma$
and $\rho=(\rho,V) \in \widehat{G^{\sigma}}$. Since $\sigma\in
Z(G^{\sigma})$, the center of $G^{\sigma}$, the Schur Lemma
implies that
\begin{equation}\label{schur}
\sigma \text{ acts by a scalar $q_{\sigma\sigma}$ on } V.
\end{equation}

\begin{lema}\label{le:odd}\cite[Lemma 2.2]{AZ}  Assume that $\sigma$ is \emph{real} (i.~e. $\sigma^{-1}\in
\oc_{\sigma}$). If $\dim\toba(\Oc_{\sigma}, \rho)< \infty$, then
$q_{\sigma\sigma} = -1$ and $s$ has even order.\qed
\end{lema}

The following result follows from \cite{H2}.

\begin{lema}\label{le:Hecke}
Let $W$ be a braided vector space, $U \subseteq W$ a braided
vector subspace of diagonal type, $\mathcal G$ the generalized
Dynkin diagram corresponding to $U$. If $\mathcal G$ contains a
$t$-cycle with $t>3$, then $\dim\toba(W)=\infty$. \qed
\end{lema}

\subsection{Symmetric groups}\label{subsec:symgroup}
Let $\sigma \in \mathbb S_m$ be of type
$(1^{n_1},2^{n_2},\dots,m^{n_m})$. Recall the notation
\eqref{eqn:descomp}. The centralizer of $\sigma$ is isomorphic to
a product $\mathbb S_m^{\sigma}=T_1 \times \cdots \times T_m, $
where
\begin{align}\label{genofcent}
T_j =\langle A_{1,j}, \dots ,A_{n_j,j} \rangle \rtimes \langle
B_{1,j}, \dots , B_{n_j-1,j} \rangle \simeq  (\Z/ j)^{n_j} \rtimes
\s_{n_j},\,
\end{align}
$1\leq j \leq m$. We will choose $A_1,\dots,A_m$ such that
$A_{1,1}=(1),\dots ,A_{n_1,1}=(n_1)$, $A_{1,2}=(n_1+1
\,\,\,n_1+2)$,\dots, $A_{n_2,2}=(n_1+2n_2-1 \,\,\,n_1+2n_2)$, and
so on. More precisely, if $1< j \leq  m$ and $r:=\sum_{1\leq k
\leq j-1}kn_k$, then
\begin{align*}
A_{l,j}&:=\Big( r+(l-1)j+1 \qquad r+(l-1)j+2 \quad \cdots \quad
r+lj \Big),
\\B_{h,j} &:= \Big( r+(h-1)j+1  \quad r+hj+1 \Big)\Big( r+(h-1)j+2
\quad r+hj+2\Big)\\ &\qquad\qquad  \cdots \Big( r+hj \quad
r+(h+1)j\Big),
\end{align*}
for all $l$, $h$, with $1\leq l \leq n_j$, $1\leq h \leq n_j-1$.
Notice that $B_{h,j}$ is an involution.

\

Let $\rho=(\rho,V) \in \widehat{\s_m^{\sigma}}$; so
\begin{align*}
\rho=\rho_1\otimes \dots \otimes \rho_m,
\end{align*}
where $\rho_j \in \widehat{T_{j}{\,}}$ has the form
\begin{equation}\label{formrho2}
\rho_j=\Ind_{\Z_{j}^{n_j} \rtimes \mathbb
S_{n_j}^{\chi_j}}^{\Z_{j}^{n_j} \rtimes \mathbb S_{n_j}} (\chi_j
\otimes \mu_j),
\end{equation} with
$\chi_j \in \widehat{\Z_{j}^{n_j}}$ and $\mu_j \in
\widehat{\s_{n_j}^{\chi_j}}$ -- see \cite[Section 8.2]{S}. Here
$\s_{n_j}^{\chi_j}$ denotes the isotropy subgroup of $\chi_j$
under the induced action of $\s_{n_j}$ over
$\widehat{\Z_{j}^{n_j}}$. Actually, $\chi_j$ is of the form
$\chi_{(t_{1,j},\dots,t_{n_j,j})},$ where $0\leq
t_{1,j},\dots,t_{n_j,j}  \leq j-1$ are such that
\begin{align}\label{lost}
\chi_{(t_{1,j},\dots,t_{n_j,j})}(A_{l,j})=\omega_j^{t_{l,j}},\qquad
1\leq l \leq n_j.
\end{align}
Notice that if $\rho_j$ is as in \eqref{formrho2}, then
\begin{align}\label{eq:degrhoj}
\deg \rho_j =[\s_{n_j}:\s_{n_j}^{(\chi_j)}] \, \deg\mu_j.
\end{align}
\smallbreak

\begin{obs}\label{obs:qeqo}
Since every $A_j$ belongs to $Z(\sm^{\sigma})$,  $A_j$ acts by a
scalar $q_{A_j}$ on $V$. Thus,
\begin{align*}
q_{\sigma\sigma}=q_{e}q_{o}, \qquad \text{where $q_{e}=\prod_{j
\text{ even}}q_{A_j}$ and  $q_{o}=\prod_{1 < j \text{
odd}}q_{A_j}$.}
\end{align*}
\end{obs}

\begin{obs}\label{subsec:1dimrep}
Assume that $\deg(\rho)=1$; that is, $\deg(\rho_j)=1$, for all
$j$. Then
\begin{align} \mathbb
S_{n_j}^{(\chi_j)}=\mathbb S_{n_j} \quad \text{ and } \quad
\mu_j=\epsilon \text{ or } \sgn \in \widehat{\mathbb S_{n_j}},
\qquad \text{ for all $j$},
\end{align}
by \eqref{eq:degrhoj}. Hence, we have that
$t_{j}:=t_{1,j}=\cdots=t_{n_j,j}$, for every $j$, and $\rho_j=
\chi_j \otimes \mu_j$. In that case, we will denote
$\chi_j=\chi_{(t_j,\dots,t_j)}$ by $\overrightarrow{\chi_{t_j}}$.
Thus, for every $j$ there exists $t_j$, with $0\leq t_j \leq j-1$,
such that
\begin{align}\label{deg=1}
\rho=(\overrightarrow{\chi_{t_1}}\otimes \mu_1) \otimes \cdots
\otimes (\overrightarrow{\chi_{t_m}}\otimes \mu_m).
\end{align}
If $n_j=0$ or $1$, then $\rho_j$ is just the trivial
representation. We will denote
\begin{align}\label{t}
{\bf t}:=(t_1,\dots,t_m),
\end{align}
which is a $m$-tuple that depends on $\rho$. Any one-dimensional
representation of $\s_m^{\sigma}$ is completely determined by
$(\mu_1,\dots,\mu_m)$ and ${\bf t}$ as above.
\end{obs}

\section{New restrictions on orbits and characters}\label{se:sobre:los:t }

Our first new result restricts the possibilities for the exponents
$t_{-,-}$'s of the representation $\rho$. The proof is an
application of the technique of abelian subracks.

\begin{lema}\label{lema:odd:gral}
Let $\rho=(\rho,V)\in \widehat{\sm^{\sigma}}$. If there exist $j$,
$l$, with $1\leq j \leq m$ and $1\leq l \leq n_j$, such that
$\omega_{j}^{4t_{l,j}}\neq 1$, then $\dim\toba(\ocs,\rho)=\infty$.
\end{lema}
\pf

Notice that $j\neq 1$, $2$. Let $N=\sum_{j\geq 3} n_j$. We
consider two cases.

\medbreak

(a) Assume that $N=1$. In this case, the type of $\sigma$ is
$(1^{n_1},2^{n_2},j)$. Then $\rho_j=\chi_{t_j}$, for some $t_j$,
$0< t_j \leq j-1$, and $q_{\sigma\sigma} =\pm \,
\omega_j^{t_j}\neq \pm 1$, by hypothesis. Now the result follows
from Lemma \ref{le:odd}.

\medbreak

(b) Assume that $N>1$. By Lemma \ref{le:odd}, we may suppose that
$q_{\sigma\sigma}= -1$. There exists $v \in V-0$ such that
$\rho(A_{l,j})v= \omega_j^{t_{l,j}} v$. We define
$\sigma_1:=\sigma$, $\sigma_2:=\sigma \,\,A_{l,j}^{-2}$,
$\sigma_3:=\sigma_2^{-1}$ and $\sigma_4:=\sigma^{-1}$; clearly,
these are four different elements. Let $\tau=(i_1 \,\,  i_2\cdots
i_j)$ be a $j$-cycle. We define
\begin{align}\label{g-1}
g_{\tau}:=\begin{cases} (i_2 \,\, i_{j}) (i_3\,\, i_{j-1}) \cdots (i_{l}\,\, i_{l+2}) &\text{, if  $j = 2l$  is even}, \\
 (i_2 \,\, i_{j}) (i_{3}\,\, i_{j-1}) \cdots (i_{l+1}\,\, i_{l+2}) & \text{, if $j = 2l + 1$ is odd.}
\end{cases}
\end{align}
Thus, $g_{\tau}$ is an involution such that $\tau^{-1}=g_{\tau}
\tau g_{\tau}$.

We choose $g_1:=\id$, $g_2:=g_{A_{l,j}}$, see \eqref{g-1},
$g_4:=g_2 g_3$  and
\begin{equation*}
g_3:=\prod_{\substack{k\neq j \\1 \leq h\leq n_k}}g_{A_{h,k}} \,
\cdot \, \prod_{\substack{1 \leq h\leq n_j\\ h\neq l}}g_{A_{h,j}}.
\end{equation*}
Then $\sigma_r=g_r \sigma g_r^{-1}$, $r=1$, $2$, $3$, $4$, and we
have the following relations
\begin{align*}
\sigma_1g_1&=g_1 \,\,\sigma_1 ,\,\,\, &\sigma_1g_2&=g_2 \,\,
\sigma_2, \,\,\, &\sigma_1g_3&=g_3 \,\, \sigma_3, \,\,\, &\sigma_1g_4&=g_4 \,\, \sigma_4,\\
\sigma_2g_1&=g_1 \,\,\sigma_2, \,\,\, &\sigma_2g_2&=g_2
\,\,\sigma_1, \,\,\, &\sigma_2g_3&=g_3 \,\,\sigma_4, \,\,\,
&\sigma_2g_4&=g_4 \,\,\sigma_3,\\
\sigma_3g_1&=g_1 \,\,\sigma_3, \,\,\, &\sigma_3g_2&=g_2
\,\,\sigma_4, \,\,\, &\sigma_3g_3&=g_3 \,\,\sigma_1, \,\,\,
&\sigma_3g_4&=g_4 \,\,\sigma_2,\\
\sigma_4g_1&=g_1 \,\,\sigma_4, \,\,\, &\sigma_4g_2&=g_2 \,\,
\sigma_3,\,\,\, &\sigma_4g_3&=g_3 \,\, \sigma_2,\,\,\,
&\sigma_4g_4&=g_4 \,\,\sigma_1.
\end{align*}
It is not difficult to see that
$W:=\ku$-span$\{g_1v,g_2v,g_3v,g_4v\}$ is a braided vector
subspace of diagonal type of $M(\mathcal O_{\sigma},\rho)$, with
braiding matrix

\begin{align*}
\mathcal Q =\begin{pmatrix}
-1 & \omega_{j}^{2t_{l,j}} & \omega_{j}^{-2t_{l,j}} & -1\\
\omega_{j}^{2t_{l,j}}  & -1 & -1 & \omega_{j}^{-2t_{l,j}}  \\
\omega_{j}^{-2t_{l,j}}  & -1 & -1 & \omega_{j}^{2t_{l,j}} \\
-1 & \omega_{j}^{-2t_{l,j}} & \omega_{j}^{2t_{l,j}} & -1
\end{pmatrix}.
\end{align*}
Since $\omega_{j}^{4t_{l,j}}\neq 1$ the generalized Dynkin diagram
is of the form given by Figure \ref{fi:rombo}. \noindent
Therefore, $\dim \toba(\mathcal O_{\sigma},\rho)=\infty$, by Lemma
\ref{le:Hecke}. \epf

\begin{figure}[ht]
\vspace{1cm}
\begin{align*}
\setlength{\unitlength}{1.4cm}
\begin{picture}(1,0)
\put(0,0){\circle*{.15}} \put(1,1){\circle*{.15}}
\put(2,0){\circle*{.15}} \put(1,-1){\circle*{.15}}
\put(0,0){\line(1,1){1}} \put(0,0){\line(1,-1){1}}
\put(1,1){\line(1,-1){1}} \put(1,-1){\line(1,1){1}}
\put(-0.2,.7){$\omega_{j}^{-4t_{l,j}}$}
\put(1.6,.7){$\omega_{j}^{4t_{l,j}}$}
\put(-0.2,-0.7){$\omega_{j}^{4t_{l,j}}$}
\put(1.6,-0.7){$\omega_{j}^{-4t_{l,j}}$}
\put(-0.4,-0.07){\small{-1}} \put(2.2,-0.07){\small{-1}}
\put(0.9,1.2){\small{-1}} \put(0.9,-1.3){\small{-1}}
\end{picture}\qquad \qquad
\end{align*}
\vspace{0.7cm}\caption{}\label{fi:rombo}
\end{figure}

\begin{obs}\label{exlema:odd}
The previous Lemma implies that if $\dim\toba(\ocs,\rho)<\infty$,
with $\rho\in \widehat{\sm^{\sigma}}$, then the scalars $q_e$ and
$q_o$ given in Remark \ref{obs:qeqo} must be $q_o=1$ and $q_e=-1$;
moreover, $t_{-,j}=0$, for all $j$ odd.
\end{obs}

Our next result discards the appearance of cycles of length $j>4$,
where $j$ is a power of 2. The proof is an application of the
technique of transversal abelian subspaces. We also need Lemma
\ref{lema:odd:gral}.

\begin{prop}\label{te:2^k}
Let $\sigma\in \s_m$, $\oc$ the conjugacy class of $\sigma$ and
$\rho \in \widehat{\s_m^{\sigma}}$. If the type of $\sigma$ is
$(1^{n_1},2^{n_2},4^{n_4},8^{n_8},\dots,(2^k)^{n_{2^k}}, \so)$,
with $k\geq 3$ and $n_{2^k}\geq 1$, then $\dim
\toba(\oc,\rho)=\infty$.
\end{prop}

\pf We may assume that $\sigma$ is of type
$(2^{n_2},4^{n_4},8^{n_8},\dots,(2^k)^{n_{2^k}})$, by \cite[Prop.
2.6]{AZ}. If $q_{\sigma\sigma}\neq -1$ or if $n_{2^k}\geq 3$, then
the result follows from Lemma \ref{le:odd} and \cite[Ex.
3.10]{AF3}, respectively. Assume that $q_{\sigma\sigma}= -1$ and
$n_{2^k}\leq 2$. We consider two cases.

(I) Assume that $n_{2^k}= 1$. Let $\alpha=(i_1\, i_2\, \cdots \,
i_{2^k})$ be the $2^k$-cycle appearing in the decomposition of
$\sigma$ as product of disjoint cycles, and we call
\begin{align*}
\mathbf I:=(i_1\, i_3\,i_5\, \cdots \, i_{2^k-1}) \quad \text{ and
} \quad {\bf P}:=(i_2\, i_4\,i_6\, \cdots \, i_{2^k}).
\end{align*}
In the proof of \cite[Lemma 2.11]{AF3}, it was shown that
\begin{itemize}
\item[(a)] ${\bf I}$ and ${\bf P}$ are disjoint $2^{k-1}$-cycles,
\item[(b)] $\alpha^2={\bf I} {\bf P}$,
\item[(c)] $\alpha  {\bf I} \alpha^{-1}={\bf P}$, (hence $\sigma  {\bf I} \sigma^{-1}={\bf P}$),
\item[(d)] ${\bf P}^t \alpha {\bf P}^{t}=\alpha^{2t+1}$, for all
integer $t$. 
\end{itemize}
For every $l \in \Z_4$, we call $\alpha_l={\bf P}^{2^{k-3}l} \,
\alpha \, {\bf P}^{-2^{k-3}l}$ and we define as in
\cite[(2.17)]{AF3},
\begin{align}\label{eq:sigmal}
\sigma_l:={\bf P}^{2^{k-3}l} \, \sigma \, {\bf P}^{-2^{k-3}l}.
\end{align}
Then $(\sigma_l)_{l \in \Z_4}$ is of type $\D_4$ in the sense of
\cite[Def. 2.2]{AF3}.

\begin{claim}
(i) $\alpha_2=\alpha^{2^{k-1}+1}$, (ii)
$\alpha_3=\alpha_1^{2^{k-1}+1}$.
\end{claim}

\pf Notice that ${\bf P}^{2^{k-2}}$ is an involution, since ${\bf
P}$ is a $2^{k-1}$-cycle. Then $\alpha_2={\bf P}^{2^{k-3}2} \,
\alpha \, {\bf P}^{-2^{k-3}2}={\bf P}^{2^{k-2}} \, \alpha \, {\bf
P}^{2^{k-2}}=\alpha^{2 \,2^{k-2}+1}=\alpha^{2^{k-1}+1}$, by (d)
above. Analogously,
\begin{align*}
\alpha_3&={\bf P}^{2^{k-3}3} \, \alpha \, {\bf P}^{-2^{k-3}3}=
{\bf P}^{2^{k-3}} {\bf P}^{2^{k-2}} \, \alpha \, {\bf
P}^{-2^{k-2}}{\bf P}^{-2^{k-3}}\\ &= {\bf P}^{2^{k-3}}\,
\alpha^{2^{k-1}+1} \, {\bf P}^{-2^{k-3}}= ({\bf P}^{2^{k-3}}\,
\alpha \, {\bf P}^{-2^{k-3}})^{2^{k-1}+1}= \alpha_1^{2^{k-1}+1}
\end{align*}
as desired. \epf

Notice that (i) implies that $\sigma_2=\sigma^{2^{k-1}+1}$ because
$\sigma^{2^{k-1}}=\alpha^{2^{k-1}}$. Analogously,
$\sigma_3=\sigma_1^{2^{k-1}+1}$. If we define
$\tau_l:=\sigma_l^{-1}$, for all $l$, then $(\sigma_l)_{l \in
\Z_4}\cup (\tau_l)_{l \in \Z_4}$ is of type $\D_4^{(2)}$. We
define $g:=(i_{1} \,\,\,i_{2^k})(i_{2} \,\,\,i_{2^k-1})\cdots
(i_{2^{k-1}-1} \,\,\,i_{2^{k-1}+1})$. Then $g$ is an involution in
$\s_m$ such that $g \trid \sigma=\sigma^{-1}$. We define
$g_l:={\bf P}^{2^{k-3}l}$ and
\begin{align*}
h_l:=g_l g,\qquad l \in \Z_4.
\end{align*}
Clearly, $h_l \trid \sigma=\tau_l$, $l \in \Z_4$.

\begin{claim}\label{le:relations}
We set $r:= 2^{k-3}$. We have the following multiplication tables:

\begin{center}
\begin{tabular}{r|cccccc}
  $\cdot $ & $g_0$& $g_1$ & $g_2$ & $g_3$ \\
  \hline
  $\sigma _0$ & $g_0 \, \sigma$ & $g_3 \, \sigma
\alpha^{2r}$ & $g_2 \, \sigma_2$ &
  $g_1 \, \sigma \alpha^{-2r}$  \\
  $\sigma _1$ & $g_2 \, \sigma \alpha^{-2r}$ & $g_1\,
\sigma$ & $g_0 \, \sigma \alpha^{2r}$ &
  $g_3 \, \sigma_2$  \\
  $\sigma_2$ & $g_0 \, \sigma_2$ & $g_3 \, \sigma
\alpha^{-2r}$ & $g_2 \, \sigma$ &
  $g_1 \, \sigma \alpha^{2r}$ \\
  $\sigma_3$ & $g_2 \, \sigma \alpha^{2r}$ & $g_1\,
\sigma_2$ & $g_0 \, \sigma \alpha^{-2r}$ &
  $g_3 \, \sigma$\\
    $\tau _0$ & $g_0 \, \sigma^{-1}$ & $g_3 \, \sigma^{-1}
\alpha^{2r}$ & $g_2 \, \sigma_2^{-1}$ &
  $g_1 \, \sigma^{-1} \alpha^{-2r}$  \\
  $\tau _1$ & $g_2 \, \sigma^{-1} \alpha^{-2r}$ & $g_1\,
\sigma^{-1}$ & $g_0 \, \sigma^{-1} \alpha^{2r}$ &
  $g_3 \, \sigma_2^{-1}$  \\
  $\tau_2$ & $g_0 \, \sigma_2^{-1}$ & $g_3 \,
\sigma^{-1} \alpha^{-2r}$ & $g_2 \, \sigma^{-1}$ &
  $g_1 \, \sigma^{-1} \alpha^{2r}$ \\
  $\tau_3$ & $g_2 \, \sigma^{-1} \alpha^{2r}$ & $g_1\,
\sigma_2^{-1}$ & $g_0 \, \sigma^{-1} \alpha^{-2r}$ &
  $g_3 \, \sigma^{-1}$
\end{tabular}
\end{center}

\begin{center}
\begin{tabular}{r|cccccc}
  $\cdot $ & $h_0$& $h_1$ & $h_2$ & $h_3$ \\
  \hline
  $\sigma _0$ & $h_0 \, \sigma^{-1}$ & $h_3
\,\sigma^{-1} \alpha^{-2r}$ & $h_2 \, \sigma_2^{-1}$ &
  $h_1 \, \sigma^{-1} \alpha^{2r}$  \\
  $\sigma _1$ & $h_2 \, \sigma^{-1} \alpha^{2r}$ & $h_1\, \sigma^{-1}$
  & $h_0 \, \sigma^{-1} \alpha^{-2r}$ &
  $h_3 \, \sigma_2^{-1}$  \\
  $\sigma_2$ & $h_0 \, \sigma_2^{-1}$ & $h_3\,
\sigma^{-1} \alpha^{2r}$ & $h_2 \, \sigma^{-1}$ &
  $h_1 \, \sigma^{-1} \alpha^{-2r}$ \\
  $\sigma_3$ & $h_2 \, \sigma^{-1} \alpha^{-2r}$ & $h_1 \, \sigma_2^{-1}$
  & $h_0 \, \sigma^{-1} \alpha^{2r}$ &
  $h_3 \, \sigma^{-1}$\\
   $\tau _0$ & $h_0 \, \sigma$ & $h_3
\,\sigma\, \alpha^{-2r}$ & $h_2 \, \sigma_2$ &
  $h_1 \, \sigma\, \alpha^{2r}$  \\
  $\tau _1$ & $h_2 \, \sigma \alpha^{2r}$ & $h_1\, \sigma$
  & $h_0 \, \sigma\, \alpha^{-2r}$ &
  $h_3 \, \sigma_2$  \\
  $\tau_2$ & $h_0 \, \sigma_2$ & $h_3\,
\sigma\, \alpha^{2r}$ & $h_2 \, \sigma$ &
  $h_1 \, \sigma\, \alpha^{-2r}$ \\
  $\tau_3$ & $h_2 \, \sigma\, \alpha^{-2r}$ & $h_1 \, \sigma_2$
  & $h_0 \, \sigma\, \alpha^{2r}$ &
  $h_3 \, \sigma$
\end{tabular}
\end{center}
\end{claim}

\pf The multiplications $\sigma_ig_j$ follow by straightforward
computations, using the fact that
$\alpha^{2^{k-1}}=\sigma^{2^{k-1}}$. For the multiplications
$\sigma_ih_j$, use that $h_{i\trid j}^{-1}\sigma_i h_j=(g_{i\trid
j}^{-1}\sigma_i g_j)^{-1}$, and the result follows. The rest can
be checked in an analogous way. Notice that $g_{i\trid
j}^{-1}\tau_i g_j=(h_{i\trid j}^{-1}\tau_i h_j)^{-1}$.\epf

Our assumption $q_{\sigma\sigma}=-1$ implies that
$\rho(\sigma_2)=-\Id$. On the other hand, we have that
$\rho(\alpha)=\omega_{2^k}^{t_{2^k}}$, and, by Lemma
\ref{lema:odd:gral}, we can suppose that
\begin{align}\label{eq:t2k}
t_{2^k}=0, \quad  2^{k-2},\quad 2^{k-1} \quad \text{or} \quad 3
\cdot 2^{k-2}.
\end{align}
Hence, $\rho(\sigma \alpha^{2r})=- \omega_{2^k}^{2rt_{2^k}}\Id=-
i^{t_{2^k}}\Id=\pm \Id$, with $i=\sqrt{-1}$, because of
\eqref{eq:t2k}. Also, $\rho(\sigma \alpha^{-2r})=-
\omega_{2^k}^{-2rt_{2^k}}\Id=-
i^{-t_{2^k}}\Id=\pm \Id$. Moreover, 
$$\rho(\sigma \alpha^{2r})=\rho(\sigma \alpha^{-2r})=-\lambda \Id,$$
with $\lambda:= i^{t_{2^k}}$. Thus, $\rho(\sigma^{-1}
\alpha^{2r})=\rho(\sigma^{-1} \alpha^{-2r})=-\lambda\Id$.

Let $v$, $w\in V-0$. We define $W:=\ku$-span of $\{u_l,v_l\,|\, l
\in \Z_4\}$, where
\begin{equation}\label{eq:uv}
\begin{aligned}
u_1&:=g_0v+g_2v, \quad &w_1&:=h_0w+h_2w,\\
u_2&:=g_0v-g_2v, \quad &w_2&:=h_0w-h_2w,\\
u_3&:=g_1v+g_3v, \quad &w_3&:=h_1w+h_3w,\\
u_4&:=g_1v-g_3v, \quad &w_4&:=h_1w-h_3w.
\end{aligned}
\end{equation}

By straightforward computations, we can see that $W$ is a braided
vector subspace of $M(\ocs,\rho)$ of Cartan type with matrix of
coefficients given by
$$\begin{pmatrix} Q&Q \\Q&Q   \end{pmatrix},\text{ where } \quad Q=\begin{pmatrix}
  -1 & -1 & -\lambda & \lambda  \\
  -1 & -1 & -\lambda & \lambda  \\
  -\lambda & \lambda & -1 & -1  \\
  -\lambda & \lambda & -1 & -1  \\
  \end{pmatrix},$$
and Dynkin diagram  given by
\begin{figure}[ht]
\vspace{1cm}
\begin{align}\label{2rombos}
\setlength{\unitlength}{1.4cm}
\begin{picture}(3,0)
\put(0,0){\circle*{.15}} \put(1,1){\circle*{.15}}
\put(2,0){\circle*{.15}} \put(1,-1){\circle*{.15}}
\put(3,0){\circle*{.15}} \put(4,1){\circle*{.15}}
\put(5,0){\circle*{.15}} \put(4,-1){\circle*{.15}}
\put(0,0){\line(1,1){1}} \put(0,0){\line(1,-1){1}}
\put(1,1){\line(1,-1){1}} \put(1,-1){\line(1,1){1}}
\put(3,0){\line(1,1){1}} \put(3,0){\line(1,-1){1}}
\put(4,1){\line(1,-1){1}} \put(4,-1){\line(1,1){1}}
\end{picture}\qquad \qquad \qquad \qquad .
\end{align}
\vspace{0.7cm}\caption{}\label{fi:rombo2}
\end{figure}

\noindent which is not of finite type. Therefore,
$\dim\toba(\ocs,\rho)=\infty$.

\bigbreak

(II) Assume that $n_{2^k}= 2$. Let $A_{1,2^k}=(i_1\, i_2\, \cdots
\, i_{2^k})$ and $A_{2,2^k}=(i_{2^k+1}\, i_{2^k+2}\, \cdots \,
i_{2^{k+1}})$ the two $2^k$-cycles appearing in $\sigma$, and let
${\bf I}={\bf I}_1 {\bf I}_2$ and ${\bf P}={\bf P}_1 {\bf P}_2$,
with
\begin{align*}
{\bf I}_1&:=(i_1\, i_3\, \cdots \, i_{2^k-1}), \quad &{\bf
I}_2&:=(i_{2^k+1}\, i_{2^k+3}\, \cdots \, i_{2^{k+1}-1}),\\
{\bf P}_1&:=(i_2\, i_4\, \cdots \, i_{2^k}),\quad &{\bf
P}_2&:=(i_{2^k+2}\, i_{2^k+4}\, \cdots \, i_{2^{k+1}}).
\end{align*}
Now, we take $\sigma_l$ as in \eqref{eq:sigmal} and we proceed in
an analogous way as in (I). \epf

\bigbreak We next examine the possible $\rho$'s when the type of
$\sigma$ is
\begin{align}
(2^{n_2},4^{n_{4}}),\quad \text{ with } n_2 \leq 5 \text{ and }
n_4 \leq 2.
\end{align}
Our first result for this question follows by performing the
argument given in \cite[Th. 4]{AF1}.

\begin{prop}\label{prop:degrho>1}
Let $\sigma\in \s_m$, $\oc$ the conjugacy class of $\sigma$ of
type $(2^{n_2},4^{n_{4}})$ and $\rho \in \widehat{\s_m^{\sigma}}$.
If $\deg\rho>1$, then $\dim \toba(\oc,\rho)=\infty$.\qed
\end{prop}

The next proposition is proved by the technique of non-abelian
subracks.

\begin{prop}\label{prop:n4=2}
Let $\sigma\in \sm$ of type $(2^{n_2},4^{2})$, $\oc$ the conjugacy
class of $\sigma$ and $\rho \in \widehat{\s_m^{\sigma}}$, with
$\deg(\rho)=1$. If $\dim \toba(\oc,\rho)<\infty$, then
$\rho_4=\chi_{(i,i)}\otimes \sgn$ or $\chi_{(-i,-i)}\otimes \sgn$,
where $i=\sqrt{-1}$.
\end{prop}

\pf We assume that $q_{\sigma\sigma}= -1$, by Lemma \ref{le:odd}.
Since $\deg(\rho)=1$, $\deg(\rho_2)=\deg(\rho_4)=1$. Let
$A_{1,4}=(j_1\,\,j_2\,\,j_3\,\,j_4)$ and
$A_{2,4}=(j_5\,\,j_6\,\,j_7\,\,j_8)$ the two 4-cycles that appear
in the decomposition of $\sigma$ as product of disjoint cycles. We
write $A_{4}=A_{1,4}A_{2,4}$. Since $\deg\rho=1$ we have that
$\rho(A_{1,4})=\rho(A_{2,4})=\omega_4^{t_4}$, with $0\leq t_4 \leq
3$ -- see Remark \ref{subsec:1dimrep}. Then
$\rho_4=\chi_{(\omega_4^{t_4},\omega_4^{t_4})}\otimes \mu_4$, with
$\mu_4=\epsilon$ \'o $\sgn$, i.~e. the trivial or sign
representation of $\Z_2$. Then $\rho(A_4)=\pm 1$. We consider two
cases.

\emph{CASE (I):} $q_{A_4}= 1$. Then
$\rho(A_{1,4})=\rho(A_{2,4})=1$ \'o $-1$. We define $s_1:=A_4$,
$s_6:=A_{1,4}A_{2,4}^{-1}$, $t_1:=A_{1,4}^{-1}A_{2,4}$,
$t_6:=s_1^{-1}$,
\begin{align*}
s_2&:=A_{1,4}\,(j_5\,\,j_6\,\,j_8\,\,j_7), &s_3&:=A_{1,4}\,(j_5\,\,j_7\,\,j_6\,\,j_8),\\
s_4&:=A_{1,4}\,(j_5\,\,j_7\,\,j_8\,\,j_6), &s_5&:=A_{1,4}\,(j_5\,\,j_8\,\,j_6\,\,j_7),\\
t_2&:= (j_1\,\,j_4\,\,j_3\,\,j_2)(j_5\,\,j_6\,\,j_8\,\,j_7),&t_3&:=(j_1\,\,j_4\,\,j_3\,\,j_2)(j_5\,\,j_7\,\,j_6\,\,j_8),\\
t_4&:=(j_1\,\,j_4\,\,j_3\,\,j_2)(j_5\,\,j_7\,\,j_8\,\,j_6),&t_5&:=(j_1\,\,j_4\,\,j_3\,\,j_2)(j_5\,\,j_8\,\,j_6\,\,j_7).
\end{align*}
Let $\alpha=\sigma s_1^{-1}$ and we define $\sigma_j:=s_j \alpha$,
$\tau_j=t_j \alpha$, $1\leq j \leq 6$. It is easy to see that
$(\sigma,\tau)$ is of type $\oct^{(2)}$ -- see \cite[Section
4]{AF3}. Now, $\tau_1=A_{1,4}^{-1}A_{2,4}\alpha=\sigma
A_{1,4}^{-2}$, $\sigma_6=A_{1,4}A_{2,4}^{-1}\alpha=\sigma
A_{2,4}^{-2}$. We choose $g:=(1\,\,2)(3\,\,4)$; then $g\trid
\sigma_1=\tau_1$, $g^{-1}\sigma_1 g=\tau_1=\sigma A_{1,4}^{-2}$ y
$g^{-1}\sigma_6 g=g^{-1}A_{1,4}A_{2,4}^{-1} g
\alpha=A_{1,4}^{-1}A_{2,4}^{-1} \alpha=\sigma s_1^{-2}$. This
implies that $\rho(\tau_1)=\rho(g^{-1}\sigma_1
g)=-\rho(A_{1,4}^{-2})=-1$,
$\rho(\sigma_6)=-\rho(A_{2,4}^{-2})=-1$ y $\rho(g^{-1}\sigma_6
g)=-\rho(s_1^{-2})=-1$. Thus, we are in the situation of \cite[Th.
4.11]{AF3}, and $\dim \toba(\oc,\rho)=\infty$.

\emph{CASE (II):} $q_{A_4}= -1$. Then
$\rho(A_{1,4})=\rho(A_{2,4})=i$ \'o $-i$, where $i=\sqrt{-1}$. Let
$\mu_4=\epsilon$. We define $s_1:=A_4$,
$$s_2:=(j_1\,\,j_2\,\,j_4\,\,j_3)(j_5\,\,j_6\,\,j_8\,\,j_7),\quad
s_3:=(j_1\,\,j_3\,\,j_2\,\,j_4)(j_5\,\,j_7\,\,j_6\,\,j_8),$$
$s_4:=s_2^{-1}$, $s_5:=s_3^{-1}$, $s_6:=s_1^{-1}$,
$t_1:=(j_1\,\,j_6\,\,j_3\,\,j_8)(j_2\,\,j_7\,\,j_4\,\,j_5)$,
$$t_2:=(j_1\,\,j_6\,\,j_4\,\,j_7)(j_2\,\,j_8\,\,j_3\,\,j_5),\quad
t_3:=(j_1\,\,j_7\,\,j_2\,\,j_8)(j_3\,\,j_6\,\,j_4\,\,j_5),$$
$t_4:=t_2^{-1}$, $t_5:=t_3^{-1}$ y $t_6:=t_1^{-1}$. Notice that
$t_1=A_{1,4}A_{2,4}B_{1,4}$, where
$B:=(j_1\,\,j_5)(j_2\,\,j_6)(j_3\,\,j_7)(j_4\,\,j_8)$.

Let $\alpha=\sigma s_1^{-1}$ and we define $\sigma_j:=s_j \alpha$,
$\tau_j=t_j \alpha$, $1\leq j \leq 6$. It is easy to see that
$(\sigma,\tau)$ is of type $\oct^{(2)}$. Now,
$\tau_1=A_{1,4}A_{2,4}B_{1,4}\alpha=\sigma B_{1,4}$,
$\sigma_6=\sigma A_{4}^{-2}$. We choose $g:=(2\, 6)(4\, 8)$; then
$g\trid \sigma_1=\tau_1$, $g^{-1}\sigma_1 g=\tau_1$ and
$g^{-1}\sigma_6 g=g^{-1}A_{1,4}^{-1}A_{2,4}^{-1}g \alpha=\sigma
A_{4}^{-2}B_{1,4}$. Then $\rho(\tau_1)=\rho(g^{-1}\sigma_1
g)=\rho(\sigma_6)=\rho(g^{-1}\sigma_6 g)=-1$, and $\dim
\toba(\oc,\rho)=\infty$, again by \cite[Th. 4.11]{AF3}. \epf

Our next task is to discard the type $(2^{2},4^{2}, \so)$; we do
this by the technique of transversal diagonal braided subspaces.

\begin{prop}\label{prop:n2=n4=2}
Let $\sigma\in \sm$ of type $(2^{2},4^{2}, \so)$, $\oc$ the
conjugacy class of $\sigma$ and $\rho \in
\widehat{\s_m^{\sigma}}$, with $\deg(\rho)=1$. Then $\dim
\toba(\oc,\rho)=\infty$.
\end{prop}

\pf By \cite[Prop. 2.6]{AZ}, we may assume that $\sigma\in
\s_{12}$ is of type $(2^{2},4^{2})$.

If $\rho_4\neq \chi_{(i,i)}\otimes \sgn$, $\chi_{(-i,-i)}\otimes
\sgn$, where $i=\sqrt{-1}$, then the result follows by Proposition
\ref{prop:n4=2}. Assume that $\rho_4=\chi_{(i,i)}\otimes \sgn$ or
$\chi_{(-i,-i)}\otimes \sgn$, with $i=\sqrt{-1}$. Following the
notation given in the preliminaries, we take $A_{1,2}=(1\,\, 2)$,
$A_{2,2}=(3\,\, 4)$, $B_{1,2}=(1\,\,3)(2\,\,4)$, $A_{1,4}=(5\,\,
6\,\, 7\,\, 8)$, $A_{2,4}=(9\,\, 10\,\, 11\,\, 12)$,
$B_{1,4}=(5\,\,9)(6\,\,10)(7\,\,11)(8\,\,12)$. We call
$A_2=A_{1,2}A_{2,2}$, $A_4=A_{1,4}A_{2,4}$ and
$\sigma=A_2A_4=(1\,\, 2)(3\,\, 4)(5\,\, 6\,\, 7\,\, 8)(9\,\,
10\,\, 11\,\, 12)$. We define $\sigma_0:=\sigma$,
$\sigma_1:=(1\,\, 2)(3\,\, 4)(5\,\, 9\,\, 7\,\, 11)(6\,\, 12\,\,
8\,\, 10)$, $\sigma_2:=\sigma_0^{-1}$, $\sigma_3:=\sigma_1^{-1}$,
\begin{align*}
\tau_0:=(1\,\, 3)(2\,\, 4)(5\,\, 6\,\, 7\,\, 8)(9\,\, 10\,\,
11\,\, 12), \tau_1:=(1\,\, 3)(2\,\, 4)(5\,\, 9\,\, 7\,\, 11)(6\,\,
12\,\, 8\,\, 10),
\end{align*}
$\tau_2:=\tau_0^{-1}$ and $\tau_3:=\tau_1^{-1}$. It is easy to see
that the family $(\sigma_l)_{l \in \Z_4}\cup (\tau_l)_{l \in
\Z_4}$ is of type $\D_4^{(2)}$. We choose $g_0:=\id$,
$g_1:=(6\,\,9)(8\,\,11)(10\,\,12)$, $g_2:=(6\,\,8)(10\,\,12)$,
$g_3:=(6\,\,11)(8\,\,9)(10\,\,12)$ and $h_l:=(2\,\,3)g_l$, $l\in
\Z_4$.

Let $v$, $w\in V-0$. We define $W:=\ku$-span of $\{u_l,v_l\,|\, l
\in \Z_4\}$, where $u_l$, $v_l$, are given by \eqref{eq:uv}. The
condition $q_{\sigma\sigma}=-1$, implies that
$\rho(A_2)=\rho_2(A_2)=1$, because $\rho_4=\chi_{(i,i)}\otimes
\sgn$ or $\chi_{(-i,-i)}\otimes \sgn$. This implies that
$\rho_2=\chi_{(t_2,t_2)}\ot \mu_{2}$, with $t_2=0$ or $1$, and
$\mu_2=\epsilon$ or $\sgn$. We consider two cases.

(a) $\rho_2=\chi_{(t_2,t_2)}\ot \epsilon$. By straightforward
computations, we can see that $W$ is a braided vector subspace of
$M(\ocs,\rho)$ of Cartan type with matrix of coefficients given by
$$\begin{pmatrix} Q&Q \\Q&Q   \end{pmatrix},\text{ where } \quad Q=\begin{pmatrix}
  -1 & -1 & -1 & 1  \\
  -1 & -1 & -1 & 1  \\
1 & -1 & -1 & -1  \\
1 & -1 & -1 & -1
  \end{pmatrix},$$
and Dynkin diagram given by \eqref{2rombos}, and
$\dim\toba(\ocs,\rho)=\infty$.

(b) $\rho_2=\chi_{(t_2,t_2)}\ot \sgn$. We proceed in an analogous
way. \epf

\medbreak Our next goal is to discard types with $n_4>0$ and
non-trivial $\so$. The proof relies on the technique of the
octahedral rack \cite[Section 4]{AF3}.

\begin{prop}\label{prop:n4>0yod}
Let $\sigma\in \sm$ of type $(1^{n_1},2^{n_2},4^{n_4}, \so)$, with
$n_4>0$ and $\so\neq \id$, $\oc$ the conjugacy class of $\sigma$
and $\rho \in \widehat{\s_m^{\sigma}}$. Then $\dim
\toba(\oc,\rho)=\infty$.
\end{prop}

\pf We assume that $q_{\sigma\sigma}= -1$, by Lemma \ref{le:odd}.
Hence $q_e=-1$ and $q_o=1$ -- see  Remark \ref{exlema:odd}. Notice
that $q_e=q_{A_2}q_{A_4}$. We consider two cases.

(I) Assume that $n_4=1$. Let $A_{1,4}=(j_1\,\,j_2\,\,j_3\,\,j_4)$
the 4-cycle appearing in the decomposition of $\sigma$ as product
of disjoint cycles. We call $s_1=A_{1,4}$,
$s_2=(j_1\,\,j_2\,\,j_4\,\,j_3)$,
$s_3=(j_1\,\,j_3\,\,j_2\,\,j_4)$, $s_4=s_2^{-1}$ and
$s_5=s_3^{-1}$, $s_6=s_1^{-1}$. Now, we define
$\sigma_l:=A_1A_2s_l\sigma_o$, $\tau_l:=A_1A_2s_l\sigma_o^{-1}$
$1\leq l \leq 6$. Then the family $(\sigma_l,\tau_l)_{1\leq l \leq
6}$ is of type $\oct^{(2)}$ -- see \cite[Def. 4.7]{AF3}.

We choose $g:=\prod_{k \text{ odd}}g_{A_k}$ -- see \eqref{g-1}.
Thus, $g$ is an involution in $\s_m$ such that $g \sigma_o g
=\sigma_o^{-1}$; then $g \sigma g =\tau_1$. Now, we compute
\begin{align*}
\rho(g^{-1} \sigma_1 g)&= \rho(\tau_1)=
\rho(A_1A_2s_1\sigma_o^{-1})=
\rho(A_2A_4)\rho(\sigma_o)^{-1}=q_eq_o \Id=-\Id,\\
\rho(\sigma_6)&=\rho(A_1A_2A_4^{-1}\sigma_o)=\rho((A_2A_4)^{-1})\rho(\sigma_o)=q_e^{-1}q_o \Id=-\Id,\\
\rho(g^{-1} \sigma_6
g)&=\rho(A_1A_2A_4^{-1}\sigma_o^{-1})=\rho((A_2A_4)^{-1})\rho(\sigma_o^{-1})=q_e^{-1}q_o^{-1}
\Id=-\Id.
\end{align*}
Then $\dim \toba(\oc,\rho)=\infty$, by \cite[Th. 4.11]{AF3}.

(II) Assume that $n_4=2$. Let $A_{1,4}=(j_1\,\,j_2\,\,j_3\,\,j_4)$
and $A_{2,4}=(j_5\,\,j_6\,\,j_7\,\,j_8)$ the two 4-cycles
appearing in the decomposition of $\sigma$. Now, we proceed as in
the previous case with $s_1=A_{1,4}A_{2,4}$,
$s_2=(j_1\,\,j_2\,\,j_4\,\,j_3)(j_5\,\,j_6\,\,j_8\,\,j_7)$,
$s_3=(j_1\,\,j_3\,\,j_2\,\,j_4)(j_5\,\,j_7\,\,j_6\,\,j_8)$,
$s_4=s_2^{-1}$, $s_5=s_3^{-1}$ and $s_6=s_1^{-1}$.\epf

\medbreak We finally discard most of the representations $\rho_1$
entering in $\rho$, see \eqref{formrho2}. We apply the technique
of $\D_3$, see \cite[Section 3]{AF3}, to the first proposition,
and the the technique of the octahedron, see \cite[Section
4]{AF3}, to the second proposition.

\begin{prop}\label{prop:degrho1>1:n2}
Let $\sigma\in \sm$ of type $(1^{n_1},2^{n_2},4^{n_4}, \so)$, with
$n_2>0$, $\oc$ the conjugacy class of $\sigma$ and $\rho \in
\widehat{\s_m^{\sigma}}$. If $\deg\rho_1>1$, then $\dim
\toba(\oc,\rho)=\infty$.
\end{prop}
\pf We assume that $q_{\sigma\sigma}= -1$, by Lemma \ref{le:odd}.
Since $\deg\rho_1>1$ we have that $n_1>0$; actually $n_1\geq 3$.
Let $A_{1,2}=(j_1\,\, j_2)$ be a transposition appearing in
$\sigma$. There exists a $j_3$ such that $\sigma$ fixes $j_3$
because $n_1>0$. We define $\sigma_1:=\sigma$, $\sigma_2:=(j_1\,\,
j_3)(j_1\,\, j_2)\sigma$ and $\sigma_3:=(j_2\,\, j_3)(j_1\,\,
j_2)\sigma$. We choose $g_1=\id$, $g_2=(j_2\,\, j_3)$ and
$g_3=(j_1\,\, j_3)$. Let $v$, $w$ be two linearly independent
vectors in $V_1$, the vector space affording $\rho_1$. We define
$W:=$span of $\{g_lv,g_lw\, | \, 1\leq l \leq 3\}$. Then $W$ is a
braided vector subspace of $M(\ocs,\rho)$ isomorphic to
$M(\oc_2^3,\sgn)\oplus M(\oc_2^3,\sgn)$, and $\dim\toba(W)=\infty$
-- see \cite[Th. 4.8]{AHS} or \cite[Th. 2.1]{AF3}. Therefore,
$\dim \toba(\oc,\rho)=\infty$. \epf

\begin{prop}\label{prop:degrho1>1:n4}
Let $\sigma\in \sm$ of type $(1^{n_1},2^{n_2},4^{n_4}, \so)$, with
$n_4>0$, $\oc$ the conjugacy class of $\sigma$ and $\rho \in
\widehat{\s_m^{\sigma}}$. If $\deg\rho_1>1$, then $\dim
\toba(\oc,\rho)=\infty$.
\end{prop}
\pf We assume that $q_{\sigma\sigma}= -1$, by Lemma \ref{le:odd}.
Analogously to the previous result, we can construct a braided
vector subspace of $M(\ocs,\rho)$ isomorphic to
$M(\oc_4^4,\chi_{(-1)})\oplus M(\oc_4^4,\chi_{(-1)})$. Then $\dim
\toba(\oc,\rho)=\infty$, by \cite[Th. 4.7]{AHS}. \epf

Our final reduction is about the characters in cases (iii) and
(iv) of  the main Theorem. We apply the technique of the rack
$\D_3^{(2)}$.

\begin{prop}\label{prop:n1>0}
Let $\sigma\in \sm$ of type $(1^{n_1},2^{n_2})$, with $n_1>0$,
$\oc$ the conjugacy class of $\sigma$ and $\rho \in
\widehat{\s_m^{\sigma}}$. If
\begin{itemize}
\item[(i)] $n_2=3$ and $\rho_2=\chi_{(3)}\otimes \epsilon$, or
\item[(ii)] $n_2=5$ and $\rho_2=\chi_{(5)}\otimes \epsilon$ or $\chi_{(5)}\otimes
\sgn$,
\end{itemize}
then $\dim \toba(\oc,\rho)=\infty$.
\end{prop}
\pf

(i) We set $\sigma=(1\, 2)(3\, 4)(5\, 6)$ and define
$\sigma_0:=\sigma$, $\sigma_1:=(1\, 2)(3\, 4)(5\, 7)$,
$\sigma_2:=\sigma_0\trid\sigma_1$, $\tau_0:=(1\, 3)(2\, 4)(5\,
6)$, $\tau_1:=\sigma_2\trid\tau_0$ and
$\tau_2:=\sigma_1\trid\tau_0$. Then $(\sigma_j,\tau_j)_{j\in
\Z_3}$ is a family of type $\D_3$ in $\oc_{\sigma}$. We choose
$g=(2\, 3)$. Thus, $\rho(g^{-1}\sigma
g)=\rho(\tau_0)=\rho_2(A_{3,2}B_{1,2}) \Id=-1$, see Subsection
\ref{subsec:symgroup}. Therefore, $\dim \toba(\oc,\rho)=\infty$,
by \cite[Th. 3.7]{AF3}.

(ii) We take $\sigma=(1\, 2)(3\, 4)(5\, 6)(7\, 8)(9\, 10)$ and
define $\sigma_0:=\sigma$, $\sigma_1:=(1\, 2)(3\, 4)(5\, 6)(7\,
8)(9\, 11)$, $\sigma_2:=\sigma_0\trid\sigma_1$, $$\tau_0:=(1\,
3)(2\, 4)(5\, 7)(6\, 8)(9\, 10),$$ $\tau_1:=\sigma_2\trid\tau_0$
and $\tau_2:=\sigma_1\trid\tau_0$. Then $(\sigma_j,\tau_j)_{j\in
\Z_3}$ is a family of type $\D_3$ in $\oc_{\sigma}$. We choose
$g=(2\, 3)(6\, 7)$. Thus, $\rho(g^{-1}\sigma
g)=\rho(\tau_0)=\rho_2(A_{5,2}B_{1,2}B_{3,2}) \Id=-1$, see
Subsection \ref{subsec:symgroup}. Hence, $\dim
\toba(\oc,\rho)=\infty$, by \cite[Th. 3.7]{AF3}. \epf

\end{document}